
\input amstex
\documentstyle{amsppt}

\hcorrection{19mm}

\nologo
\NoBlackBoxes


\topmatter
\title  Inscribing smooth knots with regular polygons
\endtitle
\author  Ying-Qing Wu$^1$
\endauthor
\leftheadtext{Ying-Qing Wu}
\rightheadtext{Completely tubing compressible tangles}
\address Department of Mathematics, University of Iowa, Iowa City, IA
52242
\endaddress
\email  wu\@math.uiowa.edu
\endemail

\thanks  $^1$ Partially supported by NSF grant \#DMS 9802558
\endthanks

\abstract A regular $n$-gon inscribing a knot is a sequence of $n$
points on a knot, such that the distances between adjacent points are
all the same.  It is shown that any smooth knot is inscribed by a
regular $n$-gon for any $n$.  
\endabstract

\endtopmatter

\document

\redefine\Gamma{G}
\define\proof{\demo{Proof}}
\define\endproof{\qed \enddemo}

\redefine\d{\delta}
\define\br{{\Bbb R}}
\define\r{\gamma}
\redefine\e{\epsilon}
\redefine\bdd{\partial}

\define\bft{{\bold t}}

\define\bfy{{\bold y}}
\define\bfz{{\bold z}}
\define\bfw{{\bold w}}
\define\bfe{{\bold e}}
\define\bfa{{\bold a}}

\redefine\Im{\text{Im}}

A knot $K: S^1 \to \br^3$ is said to be {\it inscribed by a regular
$n$-gon\/} if there is a set of points $x_0, ..., x_{n-1}$ lying on
$K$ in a cyclic order, such that the distances $\|x_{i-1} - x_{i}\|$
between $x_{i-1}$ and $x_{i}$ are the same for $i=1,...,n$, where
$x_{n} = x_0$.  Jon Simon asked the question of whether every smooth
knot $K$ is inscribed by a regular $n$-gon for all $n$.  There has
been quite some research activities on this and related problems.  See
[2, \S 11] and the references there.  In particular, it was shown by
Meyerson [4] and E.\ Kronheimer and P.\ Kronheimer [3] that given
any triangle there is one similar to it which inscribes a given planar
curve.  It is a very interesting open question whether any closed
planar curve is inscribed by a square [2], although this has been
proved for a very large class of curves, including all smooth or
piecewise linear curves [6].  See also [5].  

Up to rescaling we may assume that the length of $K$ is $1$.  It has
been observed by Eric Rawdon and Jonathan Simon (unpublished) that
given any smooth knot $K$, there is a number $N$ such that the
statement is true for all $n > N$.  Define an {\it $\epsilon$-chain\/}
of length $n$ to be a sequence of points $(x_0, ..., x_{n-1})$ on $K$,
lying successively along the positive orientation of $K$, such that
$\|x_i - x_{i-1} \| = \epsilon$ for $i=1,...,n-1$.  Choose $N$ large
enough so that for any $\epsilon < 2/N$ and any point $x$ on $K$ there
are exactly two points on $K$ with distance $\epsilon$ from $x$.  Then
starting from any $x_0 \in K$ one may construct an $\e$-chain of
length $n+1$ for small $\e$.  The sum of lengths of the short arcs on
$K$ between adjacent points in the chain is small when $\epsilon$ is
small, and exceeds $1$ when $\epsilon > 1/N$.  Hence by continuity
there must be some $\epsilon$ such that $x_{n} = x_0$, and the result
follows.  The proof fails when $n$ is small.

The question of whether every smooth knot admits an inscribed $n$-gon
for all $n$ has remained open for some time and no answer is known.  It
seems worth while to record a positive solution.  Actually, a little
more is true.  One can find a regular polygon with one vertex at any
prescribed point.  The proof is very elementary, although it does use
the concept of degree of maps between spheres in an essential way.
See [1] for some background.  Note that the theorem as stated is not
true if the smoothness assumption is dropped; however, it is not known
whether it is true if one is also allowed to move the base point.  See
Remark 6 and Conjecture 7 below for more details.

\proclaim{Theorem 1} Let $K: S^1 \to \br^3$ be a smooth knot, and let
$x_0 \in K$ be a fixed base point.  Then for any $n$ there is a
regular polygon of $n$ edges inscribing $K$, containing $x_0$ as a
vertex.  \endproclaim

\proof We may assume that $n\geq 3$; the result is trivial otherwise.
We use the same notation $K$ to denote the image of the map $K$.
Composing with $q: \br^1 \to S^1$, $q(t) = e^{2\pi ti}$, we have a
universal covering map $k: \br^1 \to K \subset \br^3$.  Without loss
of generality we may assume that $x_0$ is the origin of $\br^3$, and
$k(0) = x_0$.  Let $D(r)$ be a ball of radius $r$ centered at $x_0$,
chosen so that $D(r) \cap K = \r$ is a single arc on $K$.  Let $J$ be
the interval containing $0$ such that $k(J) = \r$.  Denote by
$d(t',t)$ the distance between two points on $K$ with parameters $t'$
and $t$, i.e., $d(t',t) = \|k(t') -k(t)\|$.

\proclaim{Lemma 2} There is a $D(r)$ and a positive number $\e$, such that 

(a) $\displaystyle \frac {\bdd d(t',t)}{\bdd t} > 0$ for $t>t'$ in
$J$;

(b) there exist $0=a_0 < a_1 < \cdots < a_{n-1}$ in $J$, such that
$d(a_{i-1}, a_i) = \e$ for $i=1,..., n-1$; and 

(c) if $d(a_{i},t)=\e$ for some $i=1,..., n-2$ and $t \in [0,1)$, then
$t=a_{i\pm 1}$.  \endproclaim

\proof (a) Since $K$ is smooth, $k'(0) \cdot k'(0) > 0$, so one can
choose $r$ small enough so that $k'(s)\cdot k'(t) > 0$ for all $s, t
\in J$.  Up to changing of coordinate we may assume that $t' =0$.  Put
$h(t) = d(t',t)^2 = d(0,t)^2 = k(t) \cdot k(t)$.  If the result were
not true then $h'(t) = 2 k(t) \cdot k'(t) = 0$.  Now $k(t) = \int_0^t
k'(s) \, ds$, hence
$$ k(t) \cdot k'(t) = \int_0^t k'(s) \cdot k'(t) \, ds = 0,$$
which is impossible because $k'(s) \cdot k'(t) > 0$ for all
$s,t$ in $J$.

(b) Choose $\e < r/n$.  Assume $a_{i}$ has been found for some $i<n-1$.
Let $b$ be the right endpoint of $J$.  Then $d(a_i,b) > r
- d(a_i,0) > r - i\e > \e$.  By (a) $d(a_i, t)$ is increasing for
$t>a_i$, and we have $d(a_i, a_i) = 0$.  Hence the continuity of
$d(a_i,t)$ implies that there is a unique $a_{i+1}$ in $(a_i, b)$
satisfying $d(a_i, a_{i+1}) = \e$.  

(c) If $k(t) \in D(r)$ the proof follows from that of (b).  If $k(t)
\notin D(r)$ one can show that $d(a_{i}, t) > \e$.
\endproof

We shall assume below that $\e$ and $a_i$ have been chosen as in the
lemma.  Let $\bfa$ denote the vector $(a_1, ..., a_{n-1})$.  Note that
$a_0 =0$ is not a component of this vector.

Let $p$ be a number between $a_{n-1}$ and $1$ such that $d(p,0) <
\e$.  This is possible because $\lim_{p\to 1} d(p,0) = 0$.  Note that
$d(0,t) > d(0,p)$ if $t\in (a_1,p)$.  Consider the set
$$ B = \{ (t_1,..., t_{n-1}) \in \br^{n-1}\,\, |\,\, a_1 \leq t_1 \leq
... \leq t_{n-1} \leq p\}.
$$ 
For notational convenience, we shall always write $t_0 = 0$ and $t_n =
1$.  Notice that $B$ is an $(n-1)$-dimensional simplex.

Denote by $\bdd B$ the boundary of $B$.  It consists of $n$ faces
given by $E = E_1 = \{\bft \in B \; | \; t_1 = a_1 \}$, $E_i = \{\bft
\in B \; | \; t_{i-1} = t_{i} \}$ for $2\leq i\leq n-1$, and $E_n =
\{\bft \in B \; | \; t_{n-1} = p \}$.  Note that $\bfa =
(a_1,...,a_{n-1})$ is an interior point of $E$.

Consider the line 
$$\Delta = \{\bfy \in \br^n \,\, | \,\, y_1 = \cdots =
y_n \},$$ 
called the diagonal of $\br^n$.  Denote by $P$ the plane 
$$P = \{\bfy \in \br^{n} \,\, | \,\, \sum y_i = 0 \}.$$
Note that $P$ is the orthogonal complement of $\Delta$ when $\Delta$
and $P$ are considered
as linear subspaces of $\br^{n}$.  Define maps
$$ 
\align
& \varphi: B \to \br^n, \qquad \varphi(\bft) = \bfy, \qquad y_i
= d(t_{i-1}, t_{i}) \\
& \psi_1: \br^{n} \to P, \qquad \psi_1(\bfy) = \bfz, \qquad z_i =
y_i - c, \qquad c = \sum y_i / n \\
& \psi_2: P - \{\bold 0 \} \to S^{n-2}, \qquad \psi_2(\bfz) = \bfz /
\|z\|
\endalign
$$
Note that $\psi_1$ is the orthogonal projection of $\br^n$ to $P$, and
$\psi_2$ is the standard radial projection.

Now suppose the theorem were not true.  Then $\varphi(B)$ is disjoint
from $\Delta$, hence the map $f = \psi_2 \circ \psi_1 \circ \varphi$
is a continuous map from $B$ to $S^{n-2}$, so its restriction to
$\bdd B$, denoted by $g: \bdd B \to S^{n-2}$, has degree $0$.  In the
following we will show that $g$ actually has degree $\pm 1$.  This
contradiction will then complete the proof of Theorem 1.

\proclaim{Lemma 3} Suppose $g(\bft) = g(\bfa)$ for some $\bft \in
\bdd B$.  Then $\bft = \bfa$.
\endproclaim

\proof Put $\bfy = \varphi(\bfa)$, and $\bfz = \psi_1(\bfy)$.  By the
choice of $\bfa$, we have $\bfy = (y_1,..., y_{n}) =
(\e,\e,...,\e,\d)$, where $\d = d(0,a_{n-1}) > d(0,a_1) = \e$ because
$d(0,t)$ is increasing on $[0,a_{n-1}]$.  Let $c = ((n-1)\e + \d)/n$
be the average of the coordinates of $\bfy$.  Then $\bfz = (\e-c, ...,
\e-c, \d-c)$.  Thus $\bfw = g(\bfa)$, the unit vector in the direction
of $\bfz$, also has the property that $w_1 = ... = w_{n-1}$, and $w_1
< w_{n}$.

Put $\bfy' = \varphi(\bft)$, and $\bfw' = g(\bft)$.  Since $\bfw' =
g(\bft) = g(\bfa) = \bfw$, we have $w'_1 = w'_i$ for
$i<n$, and $w'_1 < w'_n$.  This implies that $y'_1 = y'_i$ for $i<n$
and $y'_1 < y'_n$.  On the other hand, if $\bft \in E_i \subset \bdd
B$ for $1<i<n$ then $t_i = t_{i-1}$, so $y'_i = d(t_{i-1}, t_i) = 0$,
which is a contradiction because $t_1$ is always positive and hence
$y'_1 = d(t_1, 0) > 0$.  Also, if $i=n$, then $t_n = p$, so $y'_n =
d(0,p) \leq d(0,t_1) = y'_1$ because $t_1 \in [a_1, p]$, which again
contradicts the fact that $y'_1 < y'_n$.

Therefore we must have $\bft \in E$.  By definition we have $t_1 =
a_1$, so $y'_1 = d(t_0, t_1) = \e = y_1$.  The equalities $w'_1 =
... = w'_{n-1}$ now imply that $y'_i = d(t_{i-1},t_{i}) = y'_1 = \e$
for $i\leq n-1$.  By Lemma 2(c), this implies that $t_i = a_i$, and
the result follows.
\endproof

\proclaim{Lemma 4}  The point $\bfa \in \bdd B$ is a regular point
of $g$.
\endproclaim

\proof As usual, denote by $g_*: T_\bfa (E) \to T_{g(\bfa)}(S^{n-2})$
the map induced by $g$ on the tangent space of $E$ at $\bfa$.  We need
to show that the kernel of $g_*$ is trivial.  Put $\psi = \psi_2 \circ
\psi_1$.  Then $g_*^{-1}(0) = \varphi^{-1}_* \circ \psi^{-1}_*(0)$.
Since $\psi$ is a linear map, $\psi_* = \psi$.  One can verify that
$\psi_*^{-1}(0) = \psi^{-1}(0)$ is the linear space $L$ spanned by
$\varphi(\bfa)$ and $\bfe = (1,..., 1)$.  To prove the lemma, we need
only show that $\br^{n}$ is spanned by $L$ and $\Im(\varphi_*)$, the
image of $\varphi_*$.

From the definition of $\varphi$, one can see that $\varphi_*$ is
defined by the following matrix, which has $n$ rows and $n-1$
columns.  
$$
\pmatrix
c_{11}  & 0       & 0  & \cdots &  0  \\
c_{21}  & c_{22}  & 0  & \cdots &  0  \\
0  & c_{32} & c_{33}  & \cdots  &  0  \\
\vdots  & \vdots  & \vdots  & \cdots  &  \vdots  \\
0  & 0  & 0  & \cdots &  c_{n-1,n-1}  \\
0  & 0  & 0  & \cdots  &  c_{n,n-1} \\
\endpmatrix
$$
The coefficient $c_{ij}$ is the partial derivative of $\varphi_i$ with
respect to $t_j$ at $\bfa$.  The calculation follows from the fact
that $\varphi_i(\bft) = d(t_{i-1}, t_{i})$ is independent of $t_j$ with
$j\neq i, i-1$.  Recall from the definition of $\bfa$ that $c_{ii}$ is
positive for all $i$.

Since the boundary face $E$ of $B$ is on the plane $t_1 = a$ in
$\br^{n-1}$, the tangent space $T_{\bfa}(E)$ is spanned by the last
$n-2$ vectors of the standard basis of $\br^{n-1}$.  Hence
$\Im(\varphi_*)$ is generated by the last $n-2$ columns of the above
matrix.  It follows that $L + \Im(\varphi_*)$ is spanned by the
columns of the matrix
$$
M = 
\pmatrix
\e     & 1      & 0       & 0       & \cdots  & 0       \\
\e     & 1      & c_{22}  & 0       & \cdots  & 0       \\
\e     & 1      & c_{32}  & c_{33}  & \cdots  & 0       \\
\vdots & \vdots & \vdots  & \vdots  & \vdots  & \vdots  \\
\e     & 1      & 0       & 0       & \cdots  & c_{n-1,n-1}  \\
\d     & 1      & 0       & 0       & \cdots  & c_{n,n-1} \\
\endpmatrix
$$ where $\d = d(0, t_{n-1}) > \e$ by Lemma 2(a).  Subtracting $\e$
times of the second column from the first column and then expanding
along the first column, we have $\det M = (-1)^{n+1}
(\d-\e)c_{22}\cdots c_{n-1,n-1} \neq 0$.  Hence $\varphi$ is
transverse to $P$, and the result follows.  \endproof

Suppose $h: M \to N$ is a smooth map between closed oriented smooth
manifolds of the same dimension.  If $y\in N$ is a regular value of
$h$, then the degree of $h$ equals the number of points in $h^{-1}(y)$
at which $h$ is orientation preserving, subtracted by the number of
points in $h^{-1}(y)$ at which $h$ is orientation reversing.  See [1].
Using smooth approximation, we see that the above is true even if $h$
is only smooth in a neighborhood of $h^{-1}(y)$.  

The sphere $\bdd B$ is only piecewise smooth as a submanifold of $\Bbb
R^{n-1}$, but we can compose with a piecewise smooth map $\rho:
S^{n-2} \to \bdd B$ to obtain a map $h = g \circ \rho: S^{n-2} \to
S^{n-2}$ which is smooth in a neighborhood of $h^{-1}(g(\bfa))$.
Lemmas 3 and 4 show that $g(\bfa)$ is a regular value of $h$, and
$h^{-1}(g(\bfa))$ has only a single point.  Therefore we have $\deg(g)
= \deg(h) = \pm 1$.  As shown in the paragraph before Lemma 3, this
contradicts the facts that $g$ is the restriction of $f: B \to
S^{n-2}$ to $\bdd B$, completing the proof of Theorem 1.  \endproof

The assumptions in Theorem 1 can be weakened.  We have the following
generalization.

\proclaim{Theorem 5} Suppose $K: S^1 \to \Bbb R^m$ is a continuous
map.  Let $x_0 \in S^1$ be a point such that (i) $K$ has nonzero
continuous derivative in a neighborhood of $x_0$, and (ii) $K(x) \neq
K(x_0)$ for all $x \neq x_0$.  Then for any $n$ there is a regular
polygon of $n$ edges inscribing $K$, containing $x_0$ as a vertex.
\endproclaim

\proof Note that (i) and (ii) imply that $m\geq 2$.  By the
compactness of $S^1$ one can show that there is a neighborhood $V$ of
$x_0$ such that (ii) holds when $x_0$ is replaced by any $y \in V$.
Thus using (i) and (ii) we can still find a ball $D(r)$ of radius $r$
centered at $x_0$ such that $D(r) \cap K$ is a single arc.  In the
proof of Theorem 1 we used smoothness of $K$ in the proofs of Lemmas 2
and 4; however, in either case we only used the fact that $K$ has
continuous nonzero derivative in a small neighborhood of $x_0$.  The
rest of the proof of Theorem 1 applies verbatim.  \endproof

\example{Remark 6} Condition (i) in Theorem 5 cannot be removed.  An
easy example is given by a triangle $K$ with two equilateral edges
joined at $x_0$ at an angle less than $\pi/3$, in which case there is
no regular $3$-gon inscribing $K$ with $x_0$ as a vertex.  However, it
is not known if the result would still be true for non smooth knot if
it is not required that $x_0$ be a vertex of the polygon.

Condition (ii) in Theorem 5 cannot be removed either.  The projection
$p: \br^2 \to \br^1$ induces a smooth map $K: S^1 \to \br^1 \subset
\br^3$, which does not have a regular $3$-gon.  One can also find a
figure 8 curve on the plane with a single double point at $x_0$, which
is not inscribed by a regular 3-gon with $x_0$ as one of its vertex.
\endexample

\proclaim{Conjecture 7} Any (non smooth) knot $K: S^1 \hookrightarrow
\br^3$ is inscribed by a regular $n$-gon for any $n$.  \endproclaim

One might attempt to approach the knot $K$ with a sequence of smooth
maps $K_i$.  Let $G_i$ be a regular $n$-gon inscribing $K_i$.  The
limit of a convergent subsequence of $G_i$ is then an $n$-gon $G$ of
equal edge length inscribing $K$.  The only problem here is that $G$
might be degenerate in the sense that all of its vertices are at the
same point of $K$.


\Refs
\widestnumber\key{CGLS}

\ref \key 1 \by V. Guillemin and A. Pollack \book Differential
Topology \bookinfo Prentice-Hall \yr 1974 \endref

\ref \key 2 \by V. Klee and S. Wagon \book Old and new unsolved problems
in plane geometry and number theory  \bookinfo The Dolciani Mathematical
Expositions \vol 11 \publ Mathematical Association of America,
Washington, DC \yr 1991
\endref

\ref \key 3 \by E. Kronheimer and P. Kronheimer \paper The tripos
problem \jour J. London Math. Soc \vol 24 \yr 1981 \pages 182--192
\endref

\ref \key 4 \by Mark Meyerson \paper Equilateral triangles and
continuous curves \jour Fund. Math. \vol 110 \yr 1980 \pages 1--9
\endref

\ref \key 5 \by M. Nielsen and S. Wright \paper Rectangles inscribed
in symmetric continua \jour Geom. Dedicata \vol 56 \yr 1995 \pages
285--297 \endref

\ref \key 6 \by Walter Stromquist \paper
Inscribed squares and square-like quadrilaterals in closed curves
\jour Mathematika  \vol 36 \yr 1989 \pages 187--197
\endref

\endRefs
\enddocument